\documentclass[11pt]{article}
\usepackage{amssymb}
\newtheorem{theorem}{Theorem}
\newtheorem{corollary}{Corollary}
\begin{document}

\def\R{{\mathbb R}}
\def\C{{\mathbb C}}

\title {Integrability of the $n$--centre problem at high energies}
\author{Andreas Knauf\thanks{Mathematisches Institut,
Universit\"{a}t Erlangen-N\"{u}rnberg,
Bismarckstr.\ $1 \frac{1}{2}$, D--91054 Erlangen, Germany.
e-mail: knauf@mi.uni-erlangen.de}\and
Iskander A.\ Taimanov\thanks{Institute of Mathematics,
630090 Novosibirsk, Russia.
e-mail: taimanov@math.nsc.ru}}
\date{}
\maketitle

We consider the $n$--centre problem of celestial mechanics
in $d=2$ and $3$ dimensions.

In this paper we show that for generic configuration of the centres
at high energy levels
this system is completely integrable by using $C^\infty$ integrals
of the motion however it is not integrable in terms of
real analytic integrals.

The Hamiltonian function
\[\hat{H}:T^*\hat{M} \to \R \ \ , \ \ \hat{H}(\vec{p},\vec{q}) =
\frac{1}{2} {\vec{p}\,}^2+V(\vec{q}),\]
with potential
$$
V:\hat{M}\to\R \ \ \ ,\ \
V(\vec{q}) = -\sum_{k=1}^n\frac{Z_k}{\|\vec{q}-\vec{s}_k\|},
$$
on the cotangent bundle $T^*\hat{M}$ of configuration space
\[\hat{M} := \R^d\setminus {\{\vec{s}_1,\ldots,\vec{s}_n\}}\]
generates a -- in general incomplete -- flow.
Here $\vec{s}_k \in \R^d$ is location of the $k$-centre, $\vec{s}_k
\neq \vec{s}_l$ for $k \neq l$, and
$Z_k \in \R \setminus\{0\}$, $k=1,\dots,n$.
If $Z_k >0$ for all $k$, then this problem describes the motion of
a massive particle in the gravitational field of $n$ fixed centres.

In \cite{Knauf} it was, in particular, shown that this system admits
a smooth extension $(P,\omega,H)$ such that the corresponding flow is
complete.

Until recently it was known that

1)
for $n=1$ this system is integrable, with the angular momentum
for dimension $d=3$ being a real analytic constant of motion
(for $Z_1>0$ this is the Kepler problem);

2)
for $n=2$ this system is integrated by using
the elliptic prolate coordinates (this was done by Euler);

3)
for $n\geq3$ centres and $d=2$ it is showed in \cite{Bo}
that there is no analytic integral of the motion which
is non--constant on an energy shell $H^{-1}(E),\ E>0$;

4)
for $d=3$ and a collinear configuration of centres the angular momentum
w.r.t.\ that axis is an additional constant of the motion, independent
of the number $n$ of centres;

5)
for $d=3$ it was proved
that the topological entropy of the flow restricted to the set of
bounded orbits $b_E$ is positive
(\cite{Knauf} for
sufficiently large energies $E > E_{\mathrm{th}}$, \cite{BN}
for nonnegative energies $E \geq 0$)
and $h_{\mathrm{top}} = 0$ if $b_E$ is empty.
Furthermore $h_{\rm top}(E)$ vanishes for
$n=1$ and $2$, and $h_{\rm top}(E)>0$ if $n\geq3$ and all centres being
attracting
or not more than two $\vec{s}_k$ being on a line
(for collinear
configurations
with $Z_1,\ldots,Z_n<0$ one has $h_{\rm top}(E)=0$ for $E>0$).

Orbits of the flow fall into three classes: bounded, scattering, and
trapped. The subsets formed by these orbits are defined by $b,s$, and $t$
respectively.
The limits of scattering orbits are described by comparison
with the Kepler flow generated by the extension of
$$
\hat{H}_\infty:T^*(\R^d\setminus{\{0\}})\to\R, \ \ \
\hat{H}_\infty (\vec{p},\vec{q}) := \frac{1}{2}
\vec{p}^{\ 2} - \frac{Z_\infty}{\|\vec{q}\|}, \ \ \
Z_\infty = \sum_{k=1}^n Z_k.
$$
It was proved in \cite{Knauf} that
the set of trapped orbits is of measure zero and

\begin{itemize}
\item
for $d=2$ and attracting centres ($Z_k>0$);

\item
for $d=3$, arbitrary $Z_k \neq 0$ and
noncollinear configurations of centres
\end{itemize}

\noindent
there is a threshold energy
$E_{\mathrm{th}} \geq 0$ such that

\begin{itemize}
\item
for $E > E_{\mathrm{th}}$ many estimates are
proved and, in particular,
the set $b_E$ of bounded orbits is of measure zero;

\item
therefore above this threshold energy
almost every point $x$
in the phase space lies on a scattering orbit and
the following smooth functions are defined on he set formed by scattering
orbits:

a) the asymptotic limits of the momentum:
$$
\vec{p}^{\ \pm}: s \to \R^d;
$$

b) the time delay
$$
\tau: s \to \R
$$
which is the asymptotic difference between
the time spent by the orbit passing through $x$ and its Kepler limit
inside a ball of large radius. This function diverges near $b \cup t$.
\end{itemize}

We use these results and in the sequel assume that
$d=2$ and $Z_k>0$ or
$d=3$ and the configuration of the centres is noncollinear.

It appears that the asymptotyc limits of the momentum gives rise
to integrals of the motion. We have

\begin{theorem}
For any $E_1,E_2 > E_{\mathrm{th}}$ with $E_1 \leq E_2$,
there exists a constant
$C>0$ such that for any $g > 1$ the functions
$f_k^g: H^{-1}([E_1,E_2]) \to \R$ of the form
$$
f_k^g(x) :=
\left\{
\begin{array}{ccc}
p_k^+(x)\exp\left(-e^{\frac{C}{g-1} \sqrt{1+\tau^2(x)}}\right)&,
&x\in s\\
0&,&x\not\in s,
\end{array}
\right.
$$
are functionally independent of a full measure subset,
integrals of the motion and are of the Gevrey class of index $g$.
\end{theorem}

\begin{corollary}
On every submanifold $H^{-1}((E_1,E_2))$ with
$E_2 > E_1 > E_{\mathrm{th}}$
the $n$-centre problem is completely integrable.
\end{corollary}

As we use the same trick as the one used in
\cite{Butler,BT} it is obvious that a similar result concerning Gevrey
integrability of these systems can be obtained.
In particular, \cite{BT} an example of the integrable geodesic flow
with positive topological entropy on a compact real analytic
Riemannian manifold was constructed. Remark that
in the situation of Theorem 1 the restriction of the $n$-centre problem onto
the set of bounded orbits does not change the positive value of topological
entropy \cite{Knauf,BN}.
Moreover for large values of $E$ the $n$-centre is also not
analytically integrable as in the example given in \cite{BT}:

\begin{theorem}
On any level set $H^{-1}(E)$ with $E > E_{\mathrm{th}}$
the $n$-centre problem does not admit a pair of
functionally independent real analytic integrals of motion.
\end{theorem}

The obstruction to such an integrability is as follows.
Assume that there are such real analytic integrals of motion.
Let $P_E = H^{-1}(E)$ and
let $S$ be the subset of $P_E$
on which these integrals are
functionally dependent. It contains $b_E = P_E \cap b$, i.e. bounded
states of this energy. Take a generic point $x \in S$
and denote by $\gamma$
the intersection of $S$ with the unstable submanifold of the
Poincar\'e surface.
Fix some Riemannian metric on $P_E$.
By using results of \cite{Knauf} it is proved that for
$E > E_{\mathrm{th}}$ there have to exist
a vector $v_\infty$ which is tangent to $\gamma$ at $x$ and
a sequence $\{v_n\}$ of vectors tangent to
$P_E$ at $x$ such that
$$
\exp(x,v) \in \gamma, \ \ \ \lim_{n \to \infty} v_n = 0
$$
and
$$
\frac{\pi}{2} \geq
\angle(v_\infty,v_n) \geq O(r_n^{1+\alpha}), \ \ \ r_n = |v_n|,
$$
for some constant $\alpha \in (0,1)$.
However the analytic integrability on the level $P_E$ implies that
for a generic point $x \in S$ the set
$\gamma$ has to be a one-dimensional manifold.
By the Taylor decomposition, this implies that the angles has to converge
faster than $r_n^{1+\alpha}$:
$$
\angle(v_\infty,v_n) \sim O(r^2), \ \ \ r_n=|v_n|.
$$
Thus we arrive at a contradiction which implies Theorem 2.

Proofs of these theorems will be published elsewhere.

The second author (I.A.T.) was supported by RFBR (grant 03-01-00403) and
Max-Planck-Institute on Mathematics in Bonn.

\end{document}